\documentclass[11pt, reqno]{amsart}

 \usepackage{amsmath}
 \usepackage{amssymb}
\usepackage{bm}

\newcommand{\mysection}[1]{\section{#1}
      \setcounter{equation}{0}}

\newtheorem{theorem}{Theorem}[section]
\newtheorem{lemma}[theorem]{Lemma}

\newtheorem{corollary}[theorem]{Corollary} 

\theoremstyle{definition}
\newtheorem{assumption}{Assumption}[section]

\theoremstyle{remark}
\newtheorem{remark}{Remark}[section]

\newtheorem{example}{Example}[section]

\newcommand\bbeta{\text{\raise-.2ex\hbox{$\bm{\beta}$}}}
\newcommand\balpha{\text{
\hbox{$\bm{\alpha}$}}}

\newcommand\bA{\mathbb{A}}
\newcommand\bR{\mathbb{R}}

\newcommand\bB{\mathbb{B}}

\newcommand\bS{\mathbb{S}}

\newcommand\cB{\mathcal{B}}

\newcommand\cF{\mathcal{F}}

\newcommand\frA{\mathfrak{A}}
\newcommand\frB{\mathfrak{B}}

\newcommand\infsup{\operatornamewithlimits{inf\,\,\,sup}}
\newcommand\supinf{\operatornamewithlimits{sup\,\,\,inf}}

\begin{document}

\title[Markov policies]
{On the adjoint Markov policies in stochastic differential games}

\author{N.V. Krylov}
\email{nkrylov@umn.edu}
\address{127 Vincent Hall, University of Minnesota,
 Minneapolis, MN, 55455}

\keywords{Stochastic differential games,
Isaacs equation, value functions}

\subjclass[2010]{91A05, 91A15, 91A25}

\begin{abstract}
We consider time-homogeneous uniformly nondegenerate
stochastic differential games in domains
 and propose constructing
$\varepsilon$-optimal strategies and policies
by using adjoint Markov strategies and
adjoint Markov policies which are actually time-homogeneous Markov, however, relative
not to the original process but to a couple of processes
 governed by a system consisting of 
the main original equation and of an adjoint stochastic equations of the same type as
the main one. We show how to find $\varepsilon$-optimal strategies and policies
in these classes
by using the solvability in Sobolev spaces
 of not the original Isaacs equation
but of its appropriate modification. 
We also give an example of a uniformly nondegenerate game
 where our assumptions are not satisfied and where we conjecture that there are
no not only optimal Markov but even $\varepsilon$-optimal
adjoint (time-homogeneous) Markov strategies for one of the players.
 
\end{abstract}

\maketitle

\mysection{Introduction}
                                           \label{section 3.10.7}

Let $\bR^{d}=\{x=(x^{1},...,x^{d})\}$
be a $d$-dimensional Euclidean space and $d_{1}\geq 1$
be an integer.
Assume that we are given separable metric spaces
  $A$ and $B$,   and let,
for each $\alpha\in A$, $\beta\in B$, 
  the following 
  functions on $\bR^{d}$ are given: 

(i) $d\times d_{1} $  
matrix-valued $\sigma^{\alpha\beta}( x)
=\sigma(\alpha,\beta, x)=
(\sigma^{\alpha\beta}_{ij}( x))$,

(ii)
$\bR^{d}$-valued $b^{\alpha\beta}( x)=
b(\alpha,\beta, x)=
(b^{\alpha\beta}_{i }( x))$, and

(iii)
real-valued    $c^{\alpha\beta}( x)=c(\alpha,\beta, x)\geq0$,   
  $f^{\alpha\beta}( x)=f(\alpha,\beta, x)$, and  
$g(x)$. 
 
Under natural assumptions which will be specified later,
on a probability space
$(\Omega,\cF,P)$ carrying a $d_{1}$-dimensional Wiener process
$w_{t}$
one   associates with these objects and a bounded domain 
$G\subset\bR^{d}$ of class $C^{2}$
a stochastic differential
game 
with the diffusion term $\sigma^{\alpha\beta}(x)$,
  drift term $b^{\alpha\beta}(x)$, discount rate 
$c^{\alpha\beta}(x)$, running cost $f^{\alpha\beta}(x)$,
and the final cost $g(x)$ paid when the underlying process
first exits from $G$.
More precisely we consider the process defined by the equation
\begin{equation}
                                             \label{5.11.1}
x_{t}=x+\int_{0}^{t}\sigma^{\alpha_{s}
\beta_{s} }(  x_{s})\,dw_{s}
+\int_{0}^{t}b^{\alpha_{s}
\beta_{s} }(  x_{s})\,ds,
\end{equation}
where $\alpha_{\cdot}$ and $\beta_{\cdot}$ are admissible actions
of two players one of which is maximizing and the other minimizing an
expression like
$$
E\int_{0}^{\tau}f^{\alpha_{t}\beta_{t}}(x_{t})\,dt,
$$
where $\tau$ is the first-exit time of the process from $G$.
We adopt the setting almost identical to that of
\cite{FS_89} (although our set of admissible policies
of $\alpha$ and $\beta$ is,  generally, wider) and define
the order of players and their policies and strategies.
Then under very general conditions the value function
turns out to be a viscosity solution of the Isaacs
equation (see \cite{FS_89}). As in the case of controlled
diffusion processes and Bellman's equations it is natural
to use the Isaacs equation to construct $\varepsilon$-optimal
strategy of one player and $\varepsilon$-optimal policies of the
other. By using discrete time approximations
of this equation this was done in   \cite{FH_11}
and lead to the so-called almost optimal approximately
Markov time-inhomogeneous policies, whose actions at time $t$ depend on a very near past history.
Similar constructions one can find in \cite{Sw_96}.

In this article to find near optimal strategies and policies,
 we propose using adjoint Markov strategies and
adjoint Markov policies which are actually time-homogeneous Markov, however, relative
not to the original process $x_{t}$ but to a couple
$(x_{t},y_{t})$ which is given as a solution of a time-homogeneous
 system consisting of 
\eqref{5.11.1} and adjoint stochastic equations of the same type as
\eqref{5.11.1}. We show how to find $\varepsilon$-optimal strategies and policies
by using the solvability in Sobolev spaces
 of not the original Isaacs equation
but of its appropriate modification. Observe that it is unknown
if general even uniformly nondegenerate
 Isaacs equations have solutions in Sobolev spaces.
We also give an example of a uniformly nondegenerate game
 where our assumptions are not satisfied and where we conjecture that there are
no not only optimal Markov but even $\varepsilon$-optimal
adjoint (time-homogeneous) Markov strategy for one of the players.

As a point of comparison note that
in \cite{FS_89} and \cite{FH_11} the authors deal with time-inhomogeneous
possibly degenerate stochastic differential games
  on a finite time interval in the whole space. In our case we
have a uniformly nondegenerate
 {\em time-homogeneous\/} stochastic differential game
in a domain where it is quite natural to look
for {\em time-homogeneous\/} Markov strategies and policies.

The article is organized as follows. 
In the next section we present our main results.
In Section \ref{section 2.17.1}
we prove some auxiliary results. 
Theorems  \ref{theorem 2.10.1}
and  \ref{theorem 2.12.2}
and Lemma \ref{lemma 2.21.1} are proved in Section \ref{section 3.10.1}.
In Section \ref{section 2.23.1} we apply the previous results
to the case of controlled diffusion processes, to which
belongs Theorem \ref{theorem 2.10.2} proved in Section \ref{section 2.28.3}. 
Finally, in Section \ref{section 3.10.4} we prove Theorem \ref{theorem 3.4.1}
saying what happens if the Isaacs condition is satisfied.

By $N$ sometimes with arguments we denote various constants,
depending only on the arguments if they are present, but which
may change from one occurrence to another and, if in a statement,
we are proving, there is a claim that $N$ depends only on $a,b,...$,
then in the proof all constants called $N$ depend only on 
$a,b,...$ unless specifically indicated otherwise.

\mysection{Main results}
                                           \label{section 2.26.3}

Set $a^{\alpha\beta}=(1/2)\sigma^{\alpha\beta}\big(\sigma^{\alpha\beta}
\big)^{*}$.

\begin{assumption}
                                    \label{assumption 1.9.1}
(i) a) The functions $\sigma,b,c,f$ are continuous with respect to
$\beta\in B$ for each $(\alpha,x)$ and continuous with respect
to $\alpha\in A$ uniformly with respect to $\beta\in B$
for each $x$. b) These functions are continuous with respect to $x$
uniformly with respect to $\alpha$ and $\beta$,  the function
$g \in C^{2}(\bR^{d})$.

(ii) There are constants $K_{0}$ and $K_{1}$ such that
and for any $x,y\in\bR^{d}$ 
  $(\alpha,\beta)\in A\times B $
$$
\|\sigma^{\alpha\beta}( x)-\sigma^{\alpha\beta}( y)\|\leq K_{1}|x-y|,
\quad
  |b^{\alpha\beta}(  x)-b^{\alpha\beta}( y) |\leq K_{1}|x-y|,
$$
$$
  \|\sigma^{\alpha\beta}( x )\|,|b^{\alpha\beta}( x )|,
|c^{\alpha\beta}( x )|,|f^{\alpha\beta}( x )| 
\leq K_{0}.
$$

(iii) There is a constant $\delta\in(0,1]$
such that for any $\alpha\in A$, $\beta\in B$,  
and $x,\lambda\in\bR^{d}$ we have
$$
\delta|\lambda|^{2}\leq a^{\alpha\beta}_{ij}( x)\lambda^{i}
\lambda^{j}\leq \delta^{-1}|\lambda|^{2}.
$$
 
The reader understands, of course, that the summation
convention is adopted throughout the article.
\end{assumption}

Note that Assumption \ref{assumption 1.9.1} (iii) obviously implies that 
$d_{1}\geq d$.

Let $(\Omega,\cF,P)$ be a complete probability space,
let $\{\cF_{t},t\geq0\}$ be an increasing filtration  
of $\sigma$-fields $\cF_{t}\subset \cF$ such that
each $\cF_{t}$ is complete with respect to $\cF,P$, and let
$w_{t},t\geq0$, be a standard $d_{1}$-dimensional Wiener process
given on $\Omega$ such that $w_{t}$ is a Wiener process
relative to the filtration $\{\cF_{t},t\geq0\}$.

The following by now standard setting originated in \cite{FS_89}
although we prefer the notation introduced in \cite{Kr_13}. 
The set of progressively measurable $A$-valued
processes $\alpha_{t}=\alpha_{t}(\omega)$ is denoted by $\frA$. 
Similarly we define $\frB$
as the set of $B$-valued  progressively measurable functions.
These are the sets of policies.
By  $ \bB $ we denote
the set of (strategies) $\frB$-valued functions 
$ \bbeta(\alpha_{\cdot})$ on $\frA$
such that, for any $T\in(0,\infty)$ and any $\alpha^{1}_{\cdot},
\alpha^{2}_{\cdot}\in\frA$ satisfying
$$
P(  \alpha^{1}_{t}=\alpha^{2}_{t} 
 \quad\text{for almost all}\quad t\leq T)=1,
$$
we have
$$
P(  \bbeta_{t}(\alpha^{1}_{\cdot})=\bbeta_{t}(\alpha^{2}_{\cdot}) 
\quad\text{for almost all}\quad t\leq T)=1.
$$

For $\alpha_{\cdot}\in\frA$, 
$ \beta_{\cdot} 
\in\frB$,  and $x\in\bR^{d}$ define $x^{\alpha_{\cdot} 
\beta_{\cdot} x}_{t} $ as a unique solution of the It\^o
equation \eqref{5.11.1} and set 
$$
\phi^{\alpha_{\cdot}\beta_{\cdot} x}_{t}
=\int_{0}^{t}c^{\alpha_{s}
\beta_{s} }( x^{\alpha_{\cdot} 
\beta_{\cdot}  x}_{s})\,ds.
$$

Next, recall that $G$ is a bounded domain in $\bR^{d}$
of class $C^{2}$,
 define $\tau^{\alpha_{\cdot}\beta_{\cdot} x}$ as the first exit
time of $x^{\alpha_{\cdot} 
\beta_{\cdot} x}_{t}$ from $G$, and introduce
\begin{equation}
                                                    \label{2.12.2}
v(x)=\infsup_{\bbeta\in\bB\,\,\alpha_{\cdot}\in\frA}
E_{x}^{\alpha_{\cdot}\bbeta(\alpha_{\cdot})}\big[\int_{0}^{\tau}
f( x_{t})e^{-\phi_{t}}\,dt+g(x_{\tau})e^{-\phi_{\tau}}\big],
\end{equation}
where the indices $\alpha_{\cdot}$, $\bbeta$, and $x$
at the expectation sign are written  to mean that
they should be placed inside the expectation sign
wherever and as appropriate, that is
$$
E_{x}^{\alpha_{\cdot}\beta_{\cdot}}\big[\int_{0}^{\tau}
f( x_{t})e^{-\phi_{t}}\,dt+g(x_{\tau})e^{-\phi_{\tau}}\big]
$$
$$
:=
E \big[
 g(x^{\alpha_{\cdot}\beta_{\cdot}  x}
_{\tau^{\alpha_{\cdot}\beta_{\cdot}  x}}
)
e^{-\phi^{\alpha_{\cdot}\beta_{\cdot}  x}
_{\tau^{\alpha_{\cdot}\beta_{\cdot}  x}}}
+\int_{0}^{\tau^{\alpha_{\cdot}\beta_{\cdot}  x}}
f^{\alpha_{t}\beta_{t}   }
( 
x^{\alpha_{\cdot}\beta_{\cdot}  x}_{t})
e^{-\phi^{\alpha_{\cdot}\beta_{\cdot}  x}_{t}}\,dt\big].
$$
Observe that
this definition makes perfect sense due to 
Theorem 2.2.1  of \cite{Kr77} and
 $v(x)=g(x)$ in $\bR^{d}\setminus D$.
Similar abbreviated notation will be used in other cases
when the underlying processes and functions depend
on initial data or other parameters and functions.

Before stating our first main result   we introduce two more 
assumptions and a notation. 

\begin{assumption}
                                               \label{assumption 2.6.1}
For any $\varepsilon>0$,
there exists a finite set $\{\alpha(1),...,\alpha(n_{\varepsilon})\}\subset A$
such that for any $\alpha\in A$ there exists an $i\in\{1,...,n_{\varepsilon}\}$
such that for $u=\sigma,b,c,f$ it holds that
\begin{equation}
                                                  \label{2.6.1}
 \sup_{\substack{\beta\in B\\x\in G}}
  |u^{\alpha\beta}(x)-u^{\alpha(i)\beta}(x)|\leq \varepsilon.
\end{equation}
\end{assumption}

As is easy to see one can choose $i=i_{\varepsilon}(\alpha)$ satisfying
\eqref{2.6.1} to be a Borel function.

\begin{assumption}
                                               \label{assumption 2.12.1} 
Either $\sigma^{\alpha\beta}(x)$
are symmetric positive-definite matrix-valued functions or
 there is a constant $\nu>0$ such that $\sigma_{i,d_{1}-d+j}^{\alpha\beta}(x)=\nu
\delta_{ij}$ for all $i,j\leq d$ and all $\alpha,\beta,x$. 
\end{assumption}
The second part of 
this assumption means that the last $d$ columns of $\sigma$ form
an identity matrix multiplied by $\nu$. The only use of this 
assumption is \eqref{2.12.1} which can be satisfied in very
many other situations.

Take and fix a $\zeta\in C^{\infty}_{0}(\bR^{d})$ with unit
integral and for  a Borel measurable
$B$-valued function $\beta(\alpha,x)$ on $A\times \bR^{d}$ 
 and bounded measurable      functions $h(\alpha,\beta,x)$ given on
$A\times B\times\bR^{d}$ and $\rho>0$ set
 \begin{equation}
                                          \label{2.9.5}
h^{(\rho)}(\alpha,y,x)=\int_{\bR^{d}}h(\alpha,\beta (\alpha,y+
\rho z),x+
\rho z)
\zeta(z)\,dz,\,\,h^{(\rho)}(\alpha,y )=h^{(\rho)}(\alpha,y,y).
\end{equation}
\begin{theorem}
                                                    \label{theorem 2.10.1}
Under the above assumptions for any $\varepsilon>0$ there exist
a Borel measurable
$B$-valued function $\beta(\alpha,x)$ on $A\times \bR^{d}$
and $\rho_{0}>0$ such that, if, for $\rho\in(0,\rho_{0}]$,
 $x\in G$, and $\alpha_{\cdot}\in\frA$,
we define the process $y_{t}=y_{t}^{\alpha_{\cdot}x}(\rho)$
as a solution of
\begin{equation}
                                                            \label{2.24.01}
dy_{t}=\sigma ^{(\rho)}(\alpha_{t},y_{t} )\,dw_{t}
+b ^{(\rho)}(\alpha_{t},y_{t} )\,dt,\quad t\geq0,\quad y_{0}=x,
\end{equation}
where $\sigma^{(\rho)}$ and $b ^{(\rho)}$ are defined according to \eqref{2.9.5},
and set $\bbeta^{\rho}_{t}(\alpha_{\cdot},x)=\beta(\alpha_{t}
,y_{t}^{\alpha_{\cdot}x}(\rho))$, then
\begin{equation}
                                                              \label{2.10.08}  
v(x)\leq \sup_{\alpha_{\cdot}\in\frA}
E_{x}^{\alpha_{\cdot}\bbeta^{\rho}  (\alpha_{\cdot},x)}\big[\int_{0}^{\tau}
 f( x_{t})e^{-\phi_{t}
 }\,dt+g(x_{\tau})e^{-\phi_{\tau}
 }\big]\leq v(x)+\varepsilon.
\end{equation}

Furthermore, there exists a finite number of mutually
disjoint subsets $A_{i},i=1,...,n$, of $A$ such that $A=\bigcup_{i}A_{i}$
and for each $i$ we have   $\beta(\alpha_{1},x)=\beta(\alpha_{2},x)$ 
whenever $\alpha_{1},\alpha_{2}\in A_{i}$.
 
\end{theorem}

Observe that, obviously, \eqref{2.24.01} has a unique solution. 
Strategies like 
$$
\beta(\alpha_{t}
,y_{t}^{\alpha_{\cdot}x}(\rho))
$$
 are naturally called
adjoint Markov strategies, because their actions at time $t$ albeit are not based
only on the  current action of $\alpha$ and the current state of $x_{t}$  
 but still use  instead of the latter  the current state
of an adjoint process $y_{t}=y_{t}^{\alpha_{\cdot}x}(\rho)$, which, as we will see,
 is close to $x_{t}=x_{t}^{\alpha_{\cdot}\bbeta ^{\rho}(\alpha_{\cdot},x)x}$
if $\rho$ is small.

In the next theorem Assumption \ref{assumption 2.12.1} is not used.

\begin{theorem}
                                                    \label{theorem 2.12.2}
In Theorem \ref{theorem 2.10.1} drop Assumption \ref{assumption 2.12.1} but
suppose that on $(\Omega,\cF,P)$ there is
a Wiener process $(\hat w_{t}, \cF_{t}), t\geq0$,
independent of $w_{t}$.
 Then for any $\varepsilon>0$ there exists a constant $\nu>0$ 
such that all assertions of Theorem \ref{theorem 2.10.1} hold true
if we add to the right-hand side 
of \eqref{2.24.01} the term $\nu\, d\hat w_{t}$.

\end{theorem}

Here we see another instance of adjoint Markov strategies of the player $\beta$.
With the choice $\bbeta^{\rho}_{t}(\alpha_{\cdot},x)=\beta(\alpha_{t}
,y_{t}^{\alpha_{\cdot}x}(\rho))$ the process $x_{t}
=x_{t}^{\alpha_{\cdot}\bbeta^{\rho} (\alpha_{\cdot},x)x}$
satisfies
\begin{equation}
                                                \label{2.26.2}
dx_{t}=\sigma(\alpha_{t},\beta(\alpha_{t},y_{t}),x_{t})\,dw_{t}
+b(\alpha_{t},\beta(\alpha_{t},y_{t}),x_{t})\,dt,\quad t\geq0,\quad x_{0}=x,
\end{equation}
where $y_{t}$ is defined from \eqref{2.24.01}.
Therefore, for the player $\alpha$ to find an adequate response to
the above strategy $\bbeta^{\rho}_{t}(\alpha_{\cdot},x)$, he should
solve a more or less standard problem of optimal control of the two-component
diffusion process $(y_{t},x_{t})$ 
governed by the system \eqref{2.24.01}-\eqref{2.26.2} 
and maximize the expectation
in \eqref{2.10.08}. An unpleasant feature of this couple is that
it is always a degenerate process. It turns out that one can reduce the problem
to optimal control of only $y_{t}$ when $\rho$ is sufficiently
small and then the same Theorem \ref{theorem 2.10.1} applied in the case
of only one player will provide an adjoint Markov policy
while controlling $y_{t}$ which will become an adjoint Markov
policy of $\alpha$ in the original game. The above mentioned reduction of
the optimal control problem is based on the following.

\begin{lemma}
                                                   \label{lemma 2.21.1}
One more assertion can be added in Theorems \ref{theorem 2.10.1}
and \ref{theorem 2.12.2}: for any $\alpha_{\cdot}\in \frA$
$$
 \Big|
E_{x}^{\alpha_{\cdot}\bbeta  ^{\rho} (\alpha_{\cdot},x)}\big[\int_{0}^{\tau}
 f( x_{t})e^{-\phi_{t}}\,dt+g(x_{\tau})e^{-\phi_{\tau}}\big]
$$ 
\begin{equation}
                                                              \label{2.21.1}
-E_{x}^{\alpha_{\cdot}} \big[\int_{0}^{  \tau(\rho)}
   f ^{(\rho)}(  y_{t}(\rho))e^{-  \phi_{t}(\rho)}\,dt+g(y_{\tau(\rho)}(\rho))
e^{-  \phi_{  \tau(\rho)}(\rho)}\big]\Big|\leq \varepsilon,
\end{equation}
where  
$$
\phi^{\alpha_{\cdot}x}_{t}(\rho)=\int_{0}^{t}
  c ^{(\rho)}(\alpha_{s},y^{\alpha_{\cdot}x}_{s}(\rho))\,ds,
$$
where $f^{(\rho)}$ and $c^{(\rho)}$ are defined according to \eqref{2.9.5},
and $\tau^{\alpha_{\cdot}x}(\rho)$ is the first exit time of $y^{\alpha_{\cdot}x}_{t}
(\rho)$
from $G$.
\end{lemma}

This lemma and Theorems \ref{theorem 2.10.1} and \ref{theorem 2.12.2}
 almost immediately
lead to the following result about $\varepsilon$-optimal 
 adjoint Markov policies for $\alpha$.

\begin{theorem}
                                                    \label{theorem 2.10.2}
Let either 

(a) the assumptions of Theorem \ref{theorem 2.10.1}  be satisfied,
or 

(b) the assumptions of Theorem \ref{theorem 2.12.2} be satisfied.

Take $\varepsilon>0$, $x\in G$, $\rho$, and $\bbeta^{\rho} (\alpha_{\cdot},x)$
 from Theorem \ref{theorem 2.10.1} or \ref{theorem 2.12.2},
respectively.
Then there exist Lipschitz continuous in $x$
$d\times d_{1}$-matrix valued $\hat\sigma( x)$ and $\bR^{d}$-valued
$\hat b( x)$ given on $ \bR^{d}$, there exists a Borel measurable
$A$-valued function $\alpha^{\varepsilon}(x)$ on $ \bR^{d}$, and in case (b)
 there also exists a constant $\nu>0$,
such that, if for $x\in G$ 
we define the process $z_{t}=z^{ x}_{t}$ by
\begin{equation}
                                                              \label{2.10.7}
dz_{t}=\hat\sigma( z_{t})\,dw_{t}+
\hat b(z_{t})\,dt,\quad t\geq0,\quad z_{0}=x,
\end{equation}
in case (a) with the additional term
$\nu\, d\hat w_{t}$ on the right-hand side of \eqref{2.10.7}
in case (b)
and set $ \alpha^{\varepsilon}_{t}=\alpha^{\varepsilon}(z^{x}_{t})$, then
$$
\sup_{\alpha_{\cdot}\in\frA}
E_{x}^{\alpha_{\cdot}\bbeta ^{\rho} (\alpha_{\cdot},x)}\big[\int_{0}^{\tau}
 f( x_{t})e^{-\phi_{t}
 }\,dt+g(x_{\tau})e^{-\phi_{\tau}
 }\big]
$$
\begin{equation}
                                                                 \label{3.1.3} 
\leq  
E_{x}^{ \alpha^{\varepsilon}_{\cdot}\bbeta^{\rho} 
(  \alpha^{\varepsilon}_{\cdot},x)}\big[\int_{0}^{\tau}
 f( x_{t})e^{-\phi_{t}
 }\,dt+g(x_{\tau})e^{-\phi_{\tau}
 }\big]+\varepsilon.
\end{equation}
\end{theorem}

\begin{remark}
The above results hold  under milder assumptions
than the ones imposed. For instance, an absolutely cheep
generalization is that it suffices to have $g\in C(\bR^{d})$
rather than $g\in C^{2}(\bR^{d})$ because one can use uniform 
approximations of $g$. The domain $\Omega$ also need
not be in $C^{2}$. It is quite sufficient for it to satisfy
the exterior cone condition or be even worse than that.
Again appropriate approximations would do the job.

The point of the article was to promote adjoint Markov
policies and strategies, rather than deal
with numerous side problems arising along the way.
\end{remark}

\begin{example}
                                                  \label{example 2.16.1}
Let $d=1$, $G=(-1,1)$, $A= B=\{\pm1\}$,   $\sigma(\alpha,\beta)=\beta$,
$c=0$, $f= (1-|x+\alpha\beta|)_{+}$, $g\equiv0$. The Isaacs equation is
$$
\supinf_{\alpha\in A\,\,\beta\in B}[(1/2)u''+(1-|x+\alpha\beta|)_{+}]=0,
$$
which is equivalent to
$$
0=(1/2)u''+\supinf_{\alpha\in A\,\,\beta\in B} (1-|x+\alpha\beta|)_{+} 
=(1/2)u''.
$$
The solution of this equation in $G$ with zero boundary data is zero.
The inf inside is zero for any $\alpha$ and is
 obtained on $\beta(\alpha,x)=\alpha\,\text{sign}\,x$
($\text{sign}\,0:=-1$).

Like in \cite{FS_89} and \cite{FH_11}, let our probability space be the space $C([0,\infty))$ of
real-valued continuous functions
on $[0,\infty)$ with Wiener measure on the $\sigma$-field
of Borel subsets of $C([0,\infty))$. Let the Wiener process be defined by
$w_{t}(x_{\cdot})=x_{t}$, $t\geq0$. Also let $\cF_{t}$ 
be the $\sigma$-field generated by $w_{s},s\leq t$.

In such situation the equation
\begin{equation}
                                                                 \label{3.1.5}
dx_{t}=\text{sign}\,x_{t}\,dw_{t},\quad t\geq0,\quad x_{0}=0
\end{equation}
does not have $\cF_{t}$-adapted solutions at all (Tanaka's example),
 and $\beta$ cannot
use the strategy $\beta(\alpha,x)=\alpha\,\text{sign}\,x$, since
$\alpha$ can choose to be 1 for all times.

The author believes that in this example there is no (time-homogeneous)
$\varepsilon$-optimal adjoint Markov strategies for $\beta$
if $\varepsilon$ is small enough. Regarding time-inhomogeneous
adjoint Markov strategies the reader is referred to \cite{Kr_86}.
However, our results show that,
if we just take two independent copies of our probability space 
with $w_{t}$ being the Wiener process on one copy and 
$\hat w_{t}$ being the Wiener process on the other, take a mollification
$\chi(x)$ of $\text{sign}\,x$ take a $\nu>0$ and introduce
an adjoint process by
$$
dy_{t}=\alpha_{t}\chi(y_{t})\,dw_{t}+\nu \,d\hat w_{t},\quad t>0,\quad y_{0}=0,
$$
then the strategy $\bbeta_{t}(\alpha_{\cdot})=\alpha_{t}
\text{sign}\,y_{t}$ will be $\varepsilon$-optimal for $\beta$
if the mollification is done with kernel of sufficiently small
size and $\nu$ is sufficiently small. By the way, on
thus extended probability space \eqref{3.1.5} still does not have
solutions.
\end{example}
  
\begin{assumption}
                                                       \label{assumption 2.4.1}
 Assumption \ref{assumption 2.6.1} is not necessarily satisfied, but
for any $\varepsilon>0$,
there exists a finite set $\{\beta(1),...,\beta(n_{\varepsilon})\}\subset B$
such that for any $\beta\in B$ there exists an $i\in\{1,...,n_{\varepsilon}\}$
such that for $u=\sigma,b,c,f$ it holds that
\begin{equation}
                                                  \label{3.10.1}
 \sup_{\substack{\alpha\in A\\x\in G}}
  |u^{\alpha\beta}(x)-u^{\alpha\beta(i)}(x)|\leq \varepsilon,
\end{equation}
 and for any $u_{ij},u_{i},u$
on $G$ we have
$$
\supinf_{\alpha\in A\,\,\beta\in B}\big[a^{\alpha\beta}_{ij} u_{ij} +
b ^{\alpha\beta}_{i } u_{i} -c^{\alpha\beta} u+f^{\alpha\beta}\big]
$$
\begin{equation}
                                                             \label{3.2.2}
=\infsup_{\beta\in B\,\,\alpha\in A}\big[a^{\alpha\beta}_{ij} u_{ij} +
b ^{\alpha\beta}_{i } u_{i} -c^{\alpha\beta} u+f^{\alpha\beta}\big].
\end{equation}

\end{assumption}

When the Isaacs condition \eqref{3.2.2} is satisfied it is natural
to introduce 
   $ \bA $ as 
the set of $\frA$-valued functions 
$ \balpha(\beta_{\cdot})$ on $\frB$
such that, for any $T\in(0,\infty)$ and any $\beta^{1}_{\cdot},
\beta^{2}_{\cdot}\in\frB$ satisfying
$$
P(  \beta^{1}_{t}=\beta^{2}_{t} 
 \quad\text{for almost all}\quad t\leq T)=1,
$$
we have
$$
P(  \balpha_{t}(\beta^{1}_{\cdot})=\balpha_{t}(\beta^{2}_{\cdot}) 
\quad\text{for almost all}\quad t\leq T)=1.
$$

\begin{theorem}
                                                      \label{theorem 3.4.1}
Under the Assumptions \ref{assumption 1.9.1},
\ref{assumption 2.12.1}, and \ref{assumption 2.4.1}
for any $\varepsilon>0$ there exist
a Borel measurable
$A$-valued function $\alpha(x)$ on $ \bR^{d}$
and $\rho_{0}>0$ such that, if for $\rho\in(0,\rho_{0}]$,
 $x\in G$, and $\beta_{\cdot}\in\frB$
we define the process $y_{t}=y_{t}^{\beta_{\cdot}x}(\rho)$
as a solution of
\begin{equation}
                                                            \label{2.24.1}
dy_{t}=\sigma ^{(\rho)}(\beta_{t},y_{t} )\,dw_{t}
+b ^{(\rho)}(\beta_{t},y_{t} )\,dt,\quad t\geq0,\quad y_{0}=x,
\end{equation}
where $\sigma ^{(\rho)}$ and $b ^{(\rho)}$ are found following the example
$$
 h^{(\rho)}(\beta,y)=\int_{\bR^{d}}h(\alpha(y+\rho z),\beta,y+\rho z)
\zeta(z)\,dz,
$$
and set $\balpha^{\rho}_{t}(\beta_{\cdot},x)=\alpha( 
 y_{t}^{\beta_{\cdot}x}(\rho))$, then
\begin{equation}
                                                              \label{2.10.8}  
v(x)\geq \inf_{\beta_{\cdot}\in\frB}
E_{x}^{\balpha^{\rho}(\beta_{\cdot}x) \beta_{\cdot}  }\big[\int_{0}^{\tau}
 f( x_{t})e^{-\phi_{t}
 }\,dt+g(x_{\tau})e^{-\phi_{\tau}
 }\big]\geq v(x)-\varepsilon.
\end{equation}
\end{theorem}

\begin{remark}
Analogous theorem is valid when we 
 drop Assumption \ref{assumption 2.12.1} in Theorem \ref{theorem 3.4.1}
but
suppose that on $(\Omega,\cF,P)$ there is
a Wiener process $(\hat w_{t}, \cF_{t}), t\geq0$,
independent of $w_{t}$.

\end{remark}

\begin{remark}
 Observe that in Theorem \ref{theorem 2.10.1} we are talking
about the function $\beta(\alpha,x)$ depending both on $\alpha$ and $x$
and in Theorem \ref{theorem 3.4.1} we have a function $\alpha(x)$
of only $x$. Of course, this is because \eqref{3.10.1}
is assumed in Theorem \ref{theorem 3.4.1}. 

\end{remark}
 
\begin{remark}
 As a corollary of Theorems \ref{theorem 2.10.1} and \ref{theorem 3.4.1}
we obtain a well-known fact that our game has value and our strategies for $\beta$
and $\alpha$ form, so to speak, $\varepsilon$-saddle point
and the game may be called fair.
\end{remark}

\mysection{Auxiliary results}
                                                             \label{section 2.17.1}

Here is a well-known result which, for instance, is
 a particular case of Lemma 2.1 of \cite{Kr_13}.

\begin{lemma}
                                                   \label{lemma 2.25.2}
  Let $\sigma _{t} $  be a $d\times d_{1}$-matrix-valued  
and $b _{t} $  be an $\bR^{d}$-valued progressively measurable functions
on $\Omega\times(0,\infty) $. Suppose that
\begin{equation}
                                                \label{2.25.4}
 \|\sigma_{t} \| ,  |b_{t} |\leq K_{0},
\end{equation} 
\begin{equation}
                                                \label{2.25.3}
| \sigma_{t} ^{*}\lambda|
\geq \nu|\lambda|^{2}
\end{equation}
for all $\lambda\in\bR^{d}$ and $(\omega,t)$, where $\nu>0$
is a fixed constant. Take $x\in G$
and define $\tau$ as the first exit time from $G$ of
$$
x_{t}=x+\int_{0}^{t}\sigma _{s}\,dw_{s}+\int_{0}^{t}b_{s}\,ds.
$$
Then for any $n=1,2,...$ there exists a constant $N$,
depending only on $n$, $d$, $\nu$, $K_{0}$, and the diameter of $G$,
such that $E\tau^{n}\leq N$.

\end{lemma}

The following result is also very well known (can be obtained, for instance, 
by combining Lemma 2.8 of
\cite{GH_80} and Lemma 8.5 and Theorem 3.1 of \cite{Kr_11}).
By $\bS_{\delta}$ we denote
  the set of $d\times d$ symmetric matrices
whose eigenvalues are between $\delta$ and $\delta^{-1}$.
Introduce $D_{i}=\partial/\partial x^{i}$, $D_{ij}=D_{i}D_{j}$
and let $Du$ denote the gradient of $u$.

\begin{lemma}
                                                               \label{lemma 2.25.1}
Let $\nu\in(0,1]$. Then 
there exists a function $\Phi\in C^{2}(G)$ such that $\Phi>0$ on $G$,
$\Phi=0$ on $\partial G$, $|D\Phi|\geq 1$ on $\partial G$, and
$$
a_{ij}D_{ij}\Phi+b_{i}D_{i}\Phi\leq-1
$$
on $G$ for any $a=(a_{ij})\in \bS_{\nu}$ and $b=(b_{i})$
such that $|b|\leq K_{0}$.
\end{lemma}

The next few results are needed while investigating
how far off the adjoint processes are  of real controlled ones.

\begin{lemma}
                                                          \label{lemma 1.17.1}
Let $\sigma^{(i)}_{t}(y,x)$, $i=1,2$, be $d\times d_{1}$-matrix-valued  
and $b^{(i)}_{t}(y,x)$, $i=1,2$, be $\bR^{d}$-valued  functions
on $\Omega\times[0,\infty)\times \bR^{d}\times \bR^{d}$. 
Suppose that for each $T\in[0,\infty)$ these functions restricted
to $\Omega\times[0,T]\times \bR^{d}\times \bR^{d}$
are measurable with respect to $\cF_{T}\otimes \cB(\bR^{d})\otimes \cB(\bR^{d})$,
where $\cB(\bR^{d})$ is the Borel $\sigma$-field in $\bR^{d}$.
Assume that
$\sigma^{(i)}_{t}$ and $b^{(i)}_{t}$ are progressively measurable  for any $(x,y)$,
 $\sigma^{(1)}_{t}(y,x)$ and $b^{(1)}_{t}(y,x)$ are Lipschitz continuous
with respect to $x$ with   constant $K_{1}$,
and  $\sigma^{(2)}_{t}(y,y)$ and $b^{(2)}_{t}(y,y)$ are Lipschitz continuous
with respect to $y$ with a constant independent of $(\omega,t)$. 
Suppose that there exists a function $\Delta(y)$ on $G$ such that for any $y\in G$
\begin{equation}
                                                                 \label{2.17.4}
\|\sigma^{(1)}_{t}(y,y)-\sigma^{(2)}_{t}(y,y)\|^{2}+
|b^{(1)}_{t}(y,y)-b^{(2)}_{t}(y,y)|^{2}\leq \Delta(y)
\end{equation}
for all $(\omega,t)$. Also suppose that $\sigma_{t}^{(i)}$ and 
$b_{t}^{(i)}$ satisfy \eqref{2.25.4} and $\sigma_{t}^{(2)}$
satisfies  \eqref{2.25.3} 
for all values of indices, arguments, and all $\lambda\in\bR^{d}$.

Take $x\in G$ 
and define
the processes $x_{t}$ and $y_{t}$ by
$$
dy_{t}=\sigma^{(2)}_{t}(y_{t},y_{t})\,dw_{t}+b^{(2)}_{t}(y_{t},y_{t})\,dt,
\quad t\geq 0,\quad y_{0}=x,
$$
\begin{equation}
                                                                 \label{2.20.6}
dx_{t}=\sigma^{(1)}_{t}(y_{t},x_{t})\,dw_{t}+b^{(1)}_{t}(y_{t},x_{t})\,dt,
\quad t\geq 0,\quad x_{0}=x.
\end{equation}
Obviously this system has a unique solution. Finally, set $\theta$
to be the minimum of the exit times of $x_{t}$ and $y_{t}$ from $G$.
Then, for any    $T \in(0,\infty)$, we have
\begin{equation}
                                                                 \label{2.17.2}
E \sup_{t\leq T\wedge\theta}|x_{t}-y_{t}|^{2} \leq
 Ne^{NT} \|\Delta\|_{L_{d}(G)} ,
\end{equation}
where $N$ depends only on $d$, $\nu$, $K_{0}$,   $K_{1}$, and
the diameter of $G$.

\end{lemma}

Proof.  We modify the coefficients of system  \eqref{2.20.6} by multiplying them
by $I_{\theta> t}$, which does not affect \eqref{2.17.2}, allows us
to eliminate $\theta$ from it and also allows us to formally apply
Theorem 2.5.9 of \cite{Kr77}  
according to which the left-hand side of
\eqref{2.17.2} is less than
$$
 NTe^{NT}E\int_{0}^{T\wedge\theta}
\big(\|\sigma^{(1)}_{t}(y_{t}, y_{t})-\sigma^{(2)}_{t}( y_{t}, y_{t})\|^{2}+
|b^{(1)}_{t}( y_{t},y_{t})-b^{(2)}_{t}(y_{t},y_{t})|^{2}\big)\,dt,
$$
where $N=N(K_{1})$. In light of \eqref{2.17.4},
the expectation here is estimated by
$$
 E\int_{0}^{ \theta}
\Delta(y_{t} )\,dt 
$$
 and it only remains to apply
  Theorem 2.2.2 of \cite{Kr77}. The lemma is proved.

\begin{corollary}
                                            \label{corollary 2.17.1}
Under the assumptions of Lemma \ref{lemma 1.17.1}, For any 
$T\in(0,\infty)$,     we have
$$
 E\sup_{t\leq \theta}|x_{t}-y_{t}|^{2} \leq I+NT^{-1},
$$
where $I$ is the right-hand side of \eqref{2.17.2} and $N$
depends only on $d,\nu, K_{0}$, and the diameter of $G$.
\end{corollary}

Indeed, it suffices to use Lemma \ref{lemma 1.17.1} and observe that
$$
E\sup_{t\leq \theta}|x_{t}-y_{t}|^{2}I_{\theta>T}\leq
4\text{diam}^{2}(G)
  P 
 (\theta>T)\leq NT^{-1}E\theta\leq NT^{-1} .
$$
 
\begin{lemma}
                                                     \label{lemma 2.9.1}
Let $\sigma^{(i)}_{t}$, $b^{(i)}_{t}$, $i=1,2$, be as in Lemma
\ref{lemma 1.17.1} but independent of $(y,x)$ and assume that they
satisfy  \eqref{2.25.4} and \eqref{2.25.3} 
for all values of indices, arguments, and all $\lambda\in\bR^{d}$.    
Take
$h\in L_{d}(G)$, $x\in G$, and set 
$$
x^{(i)}_{ t}=x+\int_{0}^{t}\sigma^{(i)}_{ s}\,dw_{s}
+\int_{0}^{t}b^{(i)}_{ s}\,ds,\quad t\geq0.
$$
Introduce $\theta$ as the minimum of the first exit times
of $x^{(i)}_ {t} $, $i=1,2$,  from $G$. Let $\chi_{t}^{(i)}$, $i=1,2$,
be real-valued jointly measurable processes given on $[0,\theta]$
and bounded by a constant $K_{2}$.

 Then for any   $\kappa,\gamma>0$
$$
E\int_{0}^{\theta}|\chi_{t}^{(1)}h(x^{(1)}_{ t})-\chi_{t}^{(2)}h(x^{(2)}_{ t})|\,dt
\leq \gamma E\int_{0}^{\theta}|\chi_{t}^{(1)}-\chi_{t}^{(2)}|\,dt+
N_{1}(\gamma)+N_{2}(\kappa)
$$
\begin{equation}
                                                                 \label{2.9.2}
+N_{3}\|h\|_{L_{d}(G)} \kappa^{-2}
\big(E\sup_{t\leq \theta}|x^{(1)}_{ t} -  x^{(2)}_{ t} |^{2}\big)^{1/2},
\end{equation}
where   $N_{1}(\gamma)$ depends only on $h$, $\gamma$,
 $d$, $\nu$, $K_{0}$, and the diameter of $G$, and 
$N_{1}(\gamma)\to0$ as $\gamma\to\infty$,
 $N_{2}(\kappa)$ depends only on   $h$, $\kappa,d,\nu,K_{0}$,
and the diameter of $G$,
and $N_{2}(\kappa)\to 0$ as $\kappa\downarrow0$ and $N_{3}$
depends only on    $d$, $\nu$, $K_{0}$, and the diameter of $G$.

\end{lemma}

Proof. First observe that
$$
|\chi_{t}^{(1)}h(x^{(1)}_{ t})-\chi_{t}^{(2)}h(x^{(2)}_{ t})|\leq I+K_{2}|
h(x^{(1)}_{ t})-  h(x^{(2)}_{ t})|,
$$
where
$$
I=
|\chi_{t}^{(1)}-\chi_{t}^{(2)}|\,|h(x^{(1)}_{ t})|
\leq \gamma |\chi_{t}^{(1)}-\chi_{t}^{(2)}|+2K_{2}I_{|h(x^{(1)}_{ t})|>
\gamma}|h(x^{(1)}_{ t})|.
$$
By Theorem 2.2.2 of \cite{Kr77}
$$
E\int_{0}^{\theta}I_{|h(x^{(1)}_{ t})|>\gamma}|h(x^{(1)}_{ t})|\,dt\leq 
N\|I_{|h|>\gamma}h\|_{L_{d}(G)},
$$
where $N$ depends only on $d$, $\nu$, $K_{0}$, and the diameter of $G$.
It follows that it suffices to prove the lemma for $\chi^{(i)}\equiv 1$.

In that case we extend $h$ beyond $G$ by setting
it to be zero there, which does not affect \eqref{2.9.2},
  introduce $h^{(\kappa)}$ as the convolution of $h$
and $\kappa^{-d}\zeta(x/\kappa)$, 
and replace $h$ in the left-hand side of \eqref{2.9.2}
with $h^{(\kappa)}$. The error of the replacement is less than
$$
\sum_{i=1}^{2}E\int_{0}^{\theta}|h(x^{(i)}_{ t})-h^{(\kappa)}(x^{(i)}_{ t})|\,dt,
$$
which by Theorem 2.2.2 of \cite{Kr77} is less than a constant,
depending only on $\nu$, $d$, $K_{0}$, and the diameter of $G$, times
$$
\|h -h^{(\kappa)}\|_{L_{d}(\bR^{d})},
$$
which tends to zero as $\kappa\downarrow0$. This gives us
the   term $N_{2}(\kappa)$ on the right in \eqref{2.9.2}.
Finally,
$$
E \int_{0}^{\theta}|h^{(\kappa)}(x^{(1)}_{ t})-h^{(\kappa)}(x^{(2)}_{ t})|\,dt
\leq \sup_{\bR^{d}}|Dh^{(\kappa)}|(E\theta^{2})^{1/2}
\big(E\sup_{t\leq \theta}|x^{(1)}_{  t} -  x^{(2)}_{ t} |^{2}\big)^{1/2}
$$
$$
\leq N\kappa^{-2}\|h\|_{L_{d}(\bR^{d})}\|D\zeta\|_{L_{d/(d-1)}(\bR^{d})}
 \big(E\sup_{t\leq \theta}|x_{1,t} -  x_{2,t} |^{2}\big)^{1/2}.
$$

 The lemma is proved.

\mysection{Proof of Theorems \protect\ref{theorem 2.10.1}
and \protect\ref{theorem 2.12.2}
and Lemma \protect\ref{lemma 2.21.1}}
                                                         \label{section 3.10.1}

Recall that
 $a^{\alpha\beta}=(1/2)\sigma^{\alpha\beta}\big(\sigma^{\alpha\beta}\big)^{*}$   
and for   sufficiently smooth functions $u=u(x)$ introduce
$$
L^{\alpha\beta} u( x)=a^{\alpha\beta}_{ij}( x)D_{ij}u(x)+
b ^{\alpha\beta}_{i }( x)D_{i}u(x)-c^{\alpha\beta} ( x)u(x).
$$
 Also set
\begin{equation}
                                                     \label{1.16.1}
H[u](x)=\supinf_{\alpha\in A\,\,\beta\in B}
(L^{\alpha\beta} u(x)+f^{\alpha\beta} (x)).
\end{equation}
 
\begin{lemma}
                                               \label{lemma 1.16.1}
Take  $u  \in W^{2}_{d}(G) $ and $m\in\{1,2,...\} $.
Then for any $\alpha\in A$ there exists a Borel
 $B$-valued   function  
$\beta ( x)$  on $\bR^{d}$ such that
  for almost all $x\in G$
\begin{equation}
                                                              \label{2.8.1}
\big | L^{\alpha\beta(x)}u(x)+f^{\alpha\beta(x)}(x)-
H^{\alpha}(x)\big | \leq m^{-1},
\end{equation}
where
$$
H^{\alpha}[u](x):=\inf_{\beta\in B}\big(
L^{\alpha\beta }u(x)+f^{\alpha\beta }(x)\big).
$$

\end{lemma}

Proof. Fix $\alpha\in A$ and $u  \in W^{2}_{d}(G) $ and choose $u$, $Du$, and $D^{2}u$
so that they are Borel functions. Then let $\{\beta(i),i=1,2,...\}$ be a countable
everywhere dense set in $B$. Since $a,b,c,f$ are continuous in $\beta$,
$$
H^{\alpha}[u](x)=\inf_{\beta(i)}\big(
L^{\alpha\beta(i) }u(x)+f^{\alpha\beta(i) }(x)\big),
$$
and for any $x\in G$ there exists $\beta(i)$ with the least $i=i(x)$
for which
$$
H^{\alpha}[u](x)\geq  
L^{\alpha\beta(i(x)) }u(x)+f^{\alpha\beta (i(x))}(x)\big)-m^{-1}.
$$
As is easy to see, $i(x)$ is a Borel function and such is $\beta(i(x))$
as well. For $x\not \in G$ set $\beta(x)=\beta_{0}$, where
$\beta_{0}$ is any element of $B$. 
Then we get a function we need and the lemma is proved.

\begin{lemma}
                                               \label{lemma 2.8.1}
Take  $u  \in W^{2}_{d}(G) $ and $m\in\{1,2,...\}$.
Then there exists a finite family
of Borel $B$-valued functions $ \{\beta(1),...,\beta(n_{m})\}$
on $\bR^{d}$ and a Borel $B$-valued function $\beta(\alpha,x)$
on $A\times \bR^{d}$ such that

(i) $\beta(\alpha,\cdot)\in \{\beta(1),...,\beta(n_{m})\}$
for any $\alpha\in A$;

(ii) for
$$
h^{\alpha}=L^{\alpha\beta(\alpha,\cdot)}u+f^{\alpha\beta(\alpha,\cdot)}-H[u]
$$
we have
\begin{equation}
                                                              \label{2.8.2}
\|\sup_{\alpha\in A}h^{\alpha}_{+}\|_{L_{d}(G)}\leq m^{-1}.
\end{equation}

\end{lemma}

Proof. Again choose $u$, $Du$, and $D^{2}u$
so that they are Borel functions and take $\{\alpha(1),...,\alpha(n_{\varepsilon})\}$
from Assumption \ref{assumption 2.6.1} for $\varepsilon=1/m$. 
Then let $\beta(i,x)$
be functions found from Lemma \ref{lemma 1.16.1} corresponding to
$\alpha(i)$, $i=1,...,n_{\varepsilon}$.
Define $i(\alpha)$ to be the first $i$ for which \eqref{2.6.1} holds.
 Finally, set 
$$
\beta(\alpha,x)=\beta(i(\alpha),x).
$$

By Assumption \ref{assumption 2.6.1}, for any $\alpha\in A$, $\beta\in B$,
and $x\in G$,
$$
L^{\alpha\beta }u(x)+f^{\alpha\beta }(x)\leq 
L^{\alpha(i(\alpha))\beta }u(x)+f^{\alpha(i(\alpha))\beta }(x)
$$
$$
+m^{-1} N(1+|u(x)|+|Du(x)|+|D^{2}u(x)|),
$$
where and below the constants denoted by $N$ depend only on $d$. By plugging in $\beta=
\beta(\alpha,x)=\beta(i(\alpha),x)$ we find that, for any $\alpha\in A$
and $x\in G$
$$
h^{\alpha}(x)\leq H^{\alpha}[u](x)-H[u](x)+m^{-1} N(1+|u(x)|+|Du(x)|+|D^{2}u(x)|)
$$
$$
\leq m^{-1} N(1+|u(x)|+|Du(x)|+|D^{2}u(x)|),
$$
where the last inequality is due to $H^{\alpha}[u]\leq H[u]$.
This yields \eqref{2.8.2} with $m^{-1}$ on the right
multiplied by $N$ times the $L_{d}(G)$-norm
of $1+|u |+|Du |+|D^{2}u |$. Obviously, this is enough
and the lemma is proved.

 Set
$$
P[u](x)=\sup_{a\in\bS_{\delta}}a_{ij}D_{ij}u(x).
$$
By Theorem 14.1.6 of \cite{Kr_18} for each $K $ the equation
$$
\max(H[u_{K}],P[u_{K}]-K)=0
$$
in $G$ (a.e.) with boundary condition $u_{K}=g \in C^{2}$
has a solution $u_{K}\in W^{2}_{p}(G)$ for any $p>1$. By following the arguments
in  Section 7 of
\cite{Kr_14}, we conclude that the $u_{K}$'s admit a representation as the value
functions in the corresponding stochastic games and by Theorem 7.1 of
\cite{Kr_14} we have
 $u_{K}\downarrow v$ uniformly on 
$\bar G$  as $K\to\infty$. Observe that (a.e.) in $G$
\begin{equation}
                                                             \label{2.10.3}
H[u_{K}]\leq0.
\end{equation}

Next, fix   $K>0$ and $m\in\{1,2,...\}$. Below we   introduce
some objects which may change as we change $K$ and $m$, but 
we still do not exhibit their dependence on $K,m$ for simplicity
of notation and because $K,m$ are fixed for now.

Let  $\{\beta (1),...,\beta (n )\}$ and
$\beta (\alpha,x)$  
be the family of functions $\beta(i)$ and function  $\beta(\alpha,x)$
from Lemma \ref{lemma 2.8.1}
with $u_{K}$   in place of $u$. Observe that by construction and \eqref{2.10.3}
\begin{equation}
                                                             \label{2.10.4}
\sup_{\alpha\in A}\big(L^{\alpha\beta (\alpha,\cdot)}u_{K}+
f^{\alpha\beta (\alpha,\cdot)}\big)\leq h,
\end{equation}
where $h\geq0$ is such that $\|h\|_{L_{d}(G)}\leq 1/m$.

    Use this $\beta(\alpha,x)$ in
\eqref{2.9.5} and \eqref{2.24.01} to define $y_{t}=y_{t}^{\alpha_{\cdot}x}(\rho)$,
$\bbeta^{\rho}_{t}(\alpha_{\cdot},x)=\beta(\alpha_{t}
,y_{t}^{\alpha_{\cdot}x}(\rho))$, and $x_{t}=x_{t}^{\alpha_{\cdot}   
\bbeta^{\rho} (\alpha_{\cdot},x)x}$. First, we want to prove that
$x_{t}$ and $y_{t}$ are close when $\rho$ is sufficiently small.
This will be based in part on the fact that the couple $(y_{t},x_{t})$
is a solution of the system 
\begin{equation}
                                                      \label{2.8.3}
\begin{aligned}
dx_{t}=& \sigma(\alpha_{t},\beta (\alpha_{t},y_{t}),x_{t})\,dw_{t}
+b(\alpha_{t},\beta (\alpha_{t},y_{t}),x_{t})\,dt,\\
dy_{t}=&\sigma ^{(\rho)}(\alpha_{t},y_{t},y_{t})\,dw_{t}
+b ^{(\rho)}(\alpha_{t},y_{t},y_{t})\,dt.
\end{aligned}
\end{equation}

An important and easy consequence of Assumption \ref{assumption 2.12.1}
is that
\begin{equation}
                                                 \label{2.12.1}
\sigma ^{(\rho)}(\alpha,y)(\sigma ^{(\rho)}(\alpha,y))^{*}\geq \nu^{2}(\delta_{ij}),
\end{equation}
for all $\rho,\alpha,y $.

\begin{lemma}
                                                        \label{lemma 2.20.2}

For any vector-valued $h=h(\alpha,\beta,x)$ define
$$
I^{h}_{\rho}(\hat \alpha ,\alpha ,y):=
  \big|h(\hat\alpha ,\beta (\alpha ,y ),y)-
\int_{\bR^{d}}h(\hat\alpha ,\beta (\alpha ,y+
\rho z),y+
\rho z)
\zeta(z)\,dz\big |^{2}.
$$
Then for any $\varepsilon>0$ there exist  $\rho_{0}>0$ and a function
$\Delta^{h}_{\rho}(y)$ such that, for all
$\rho\in(0,\rho_{0}]$,  $\hat\alpha,\alpha\in A$, $y\in G$, and
$h=\sigma,b,c,f$ we have
\begin{equation}
                                                                \label{2.20.5}
I^{h}_{\rho}(\hat \alpha ,\alpha ,y) 
\leq \Delta^{h}_{\rho}(y),\quad \|\Delta^{h}_{\rho}\|_{L_{d}(G)}\leq\varepsilon.
\end{equation}

\end{lemma}

Proof.
According to Assumption \ref{assumption 2.6.1} for any $\varepsilon>0$
there exists a finite subset $\hat A(\varepsilon)$ (independent of $\rho$) of $A$ such that
$$
\sup_{\hat\alpha\in A}I^{h}_{\rho}(\hat \alpha ,\alpha ,y)
\leq \sup_{\hat\alpha\in \hat A(\varepsilon)}I^{h}_{\rho}(\hat \alpha ,\alpha ,y)+\varepsilon
\leq\sum_{\hat\alpha\in \hat A(\varepsilon)}I^{h}_{\rho}(\hat \alpha ,\alpha ,y)+\varepsilon.
$$

 Take an $\hat\alpha\in \hat A(\varepsilon)$ and observe that the set
 $$
S(\hat \alpha):=\{ h(\hat\alpha,\beta (\alpha,\cdot),\cdot);
\alpha\in A\}
$$
 is finite (see Lemma \ref{lemma 2.8.1}) and each element of this set
is bounded and measurable
with respect to $y$. By the Lebesgue theorem
$$
I^{h}_{\rho}(\hat \alpha ,\alpha ,y)\leq  
\sum_{g\in S(\hat\alpha)}
\big|g(y)-\int_{\bR^{d}}g(y+\rho z)\zeta(z)\,dz\big |^{2}\to0
$$
as $\rho\downarrow0$ at almost any point $y\in\bR^{d}$. Hence,
$$
I^{h}_{\rho}(\hat \alpha ,\alpha ,y)\leq
\sum_{\hat\alpha\in \hat A(\varepsilon)}\sum_{g\in S(\hat\alpha)}
\big|g(y)-\int_{\bR^{d}}g(y+\rho z)\zeta(z)\,dz\big |^{2}+\varepsilon
=:\Delta^{h}_{\rho,\varepsilon}( y)+\varepsilon,
$$
where
$\Delta^{h}_{\rho,\varepsilon}$
 are bounded uniformly with respect to
$\rho$ and tend to zero as $\rho\downarrow0$ (a.e.) in $\bR^{d}$,
in particular, in $L_{d}(G)$ for any   $\varepsilon$. As a result,
for any  $\hat\alpha,\alpha\in A$ and $y\in G$,
$$
 I^{h}_{\rho}(\hat \alpha ,\alpha ,y)
\leq  
\Delta^{h}_{\rho,\varepsilon}( y)+ \varepsilon, 
$$
where for all sufficiently small $\rho$
$$
\|\Delta^{h}_{\rho,\varepsilon} + \varepsilon\|_{L_{d}(G)}
\leq 2\varepsilon N(d)\text{diam}(G).
$$
This is, certainly, enough and the lemma is proved.

\begin{lemma}
                                               \label{lemma 2.9.2}
Introduce $\theta=\theta^{\alpha_{\cdot}\bbeta ^{\rho}(\alpha_{\cdot},x)x}(\rho)$ 
 as the minimum of the first exit times
of $x_{ t}^{\alpha_{\cdot}\bbeta ^{\rho}(\alpha_{\cdot},x)x}$
and of $y_{t} ^{\alpha_{\cdot}x}(\rho)$  from $G$. Then
\begin{equation}
                                                          \label{2.18.1}
\sup_{\alpha_{\cdot}\in\frA}E_{x}^{\alpha_{\cdot}\bbeta ^{\rho}
(\alpha_{\cdot},x)}
 \sup_{t\leq \theta(\rho)}|x_{t}-y_{t}(\rho)|^{2} \to0
\end{equation}
as $\rho\downarrow0$ uniformly with respect to $x\in G$.
\end{lemma}

Proof.  By Corollary \ref{corollary 2.17.1} and Lemma \ref{lemma 2.20.2},
for any $\varepsilon,T>0$
the left-hand side of \eqref{2.18.1} is less than
$Ne^{NT}\varepsilon+N/T$, where $N$ is independent of $\rho,\varepsilon, T$,
 for all small enough $\rho$ and so is
its lim sup   as $\rho\downarrow0$.
Sending first $\varepsilon\downarrow0$ and then $T\to\infty$
yields the desired result. The lemma is proved.

\begin{corollary}
                                                 \label{corollary 2.26.1}
For $h=\sigma,b,c,f$ we have
\begin{equation}
                                                         \label{2.26.6}
\sup_{\alpha_{\cdot}\in\frA}E_{x}^{\alpha_{\cdot}\bbeta ^{\rho}
(\alpha_{\cdot},x)}
\int_{0}^{\theta}|h(\alpha_{t},\beta(\alpha_{t},y_{t}(\rho),x_{t})
-h^{(\rho)}(\alpha_{t},y_{t}(\rho),y_{t}(\rho))| \,dt\to0
\end{equation}
as $\rho\downarrow 0$ uniformly with respect to $x\in G$,
where $h^{(\rho)}(\alpha,y,x)$ is introduced according to  \eqref{2.9.5}.

\end{corollary}

Indeed, since $h$ is continuous in $x$ uniformly with respect to
$(\alpha,\beta)$, one can replace $x_{t}$ in \eqref{2.26.6}
with $y_{t}(\rho)$ only incurring the error
$$
\sup_{\alpha_{\cdot}\in\frA}E_{x}^{\alpha_{\cdot}\bbeta ^{\rho}
(\alpha_{\cdot},x)}\theta w\big(\sup_{t\leq \theta(\rho)}|x_{t}-y_{t}(\rho)|
\big)
$$ 
\begin{equation}
                                                               \label{2.27.4}
\leq \Big(\sup_{\alpha_{\cdot}\in\frA}E_{x}^{\alpha_{\cdot}\bbeta ^{\rho}
(\alpha_{\cdot},x)}\theta^{2}\Big)^{1/2}
 \Big(\sup_{\alpha_{\cdot}\in\frA}E_{x}^{\alpha_{\cdot}\bbeta ^{\rho}
(\alpha_{\cdot},x)} w^{2}\big(\sup_{t\leq \theta(\rho)}|x_{t}-y_{t}(\rho)|\Big)^{1/2},
\end{equation}
where $w(r)$, $r\geq0$, is a bounded continuous function, $w(0)=0$.
By Lemmas \ref{lemma 2.25.2} and \ref{lemma 2.9.2} this error
tends to zero as $\rho\downarrow 0$ uniformly with respect to $x\in G$.
Due to Theorem 2.2.2 of \cite{Kr77}
and Lemma \ref{lemma 2.20.2}, what remains after the above mentioned replacement
is less than a constant independent of $\rho$ times
the $L_{d}$-norm of $\Delta^{h}_{\rho}$, which also
tends to zero as $\rho\downarrow 0$ uniformly with respect to $x\in G$.

\begin{theorem}
                                                        \label{theorem 1.11.1}

For any
$x\in G$,  $\rho,\gamma,\kappa >0$   we have
$$
u_{K}(x)\geq \sup_{\alpha_{\cdot}\in\frA}
E_{x}^{\alpha_{\cdot}\bbeta ^{\rho}(\alpha_{\cdot},x)}\big[\int_{0}^{\tau}
 f( x_{t})e^{-\phi_{t}
 }\,dt+g(x_{\tau})e^{-\phi_{\tau}
 }\big]
$$
\begin{equation}
                                                             \label{2.10.6}
-\mu(\rho)(1+\gamma+\kappa^{-2}) -N_{1}(\gamma)-N_{2}(\kappa)-N m^{-1},
\end{equation}
where $N_{1}(\gamma)$ is independent of $\rho,\kappa $, 
$N_{1}(\gamma)\to0$ as $\gamma\to\infty$, $N_{2}(\kappa)$
is independent of $\rho$,
$N_{2}(\kappa)\to 0$ as $\kappa\downarrow0$, 
$N $ depends only on $d,\delta,K_{0}$, and the diameter of $G$,
$\mu(\rho)$ is independent of $\gamma,\kappa$
and $\mu(\rho)\to0$ as $\rho\downarrow0$.
\end{theorem}

Proof. For simplicity of notation
we drop the argument $\rho$ of $\theta$ and $y_{t}$.
Take $\alpha_{\cdot}\in\frA$ and observe that in the 
notation from Lemma \ref{lemma 2.9.2}
by It\^o's formula
$$
E_{x}^{\alpha_{\cdot}\bbeta ^{\rho}(\alpha_{\cdot},x)}
u_{K}(x_{\theta })e^{-\phi_{\theta }}=u_{K}(x)
$$
$$
+
E_{x}^{\alpha_{\cdot}\bbeta ^{\rho}(\alpha_{\cdot},x)}\int_{0}^{\theta}
\big[a^{ij} (\alpha_{t},\beta (\alpha_{t},y_{t} ),x_{t})D_{ij}u_{K}(x_{t})
$$
\begin{equation}
                                                              \label{2.10.1}
+b^{i}(\alpha_{t},\beta (\alpha_{t},y_{t}),x_{t})D_{i}u_{K}(x_{t})
-c(\alpha_{t},\beta (\alpha_{t},y_{t}),x_{t}) u_{K}(x_{t}) \big]e^{-\phi_{t}}\,dt,
\end{equation}
where, dropping obvious values of indices,
$$
\phi_{t}=\int_{0}^{t}c(\alpha_{s},\beta (\alpha_{s},y_{s}),x_{s})\,ds.
$$

By Lemma \ref{lemma 2.9.1} with $h=u_{K},Du_{K},D^{2}u_{K}$, for any $\kappa,\gamma>0$,
 the last term
in \eqref{2.10.1} is less than
$$
E_{x}^{\alpha_{\cdot} }\int_{0}^{\theta}
\big[a^{ij} (\alpha_{t},\beta (\alpha_{t},y_{t}),y_{t})D_{ij}u_{K}(y_{t})
$$
$$
+b^{i}(\alpha_{t},\beta (\alpha_{t},y_{t}),y_{t})D_{i}u_{K}(y_{t})
-c(\alpha_{t},\beta (\alpha_{t},y_{t}),y_{t}) u_{K}(y_{t}) \big]e^{-\phi_{t}}\,dt
$$
$$
+\gamma [I^{a}(\alpha_{\cdot},\rho, x)+
I^{b}(\alpha_{\cdot},\rho, x)+I^{c}(\alpha_{\cdot},\rho, x)]
$$
\begin{equation}
                                                              \label{2.10.01}
+N_{1}(\gamma)+N_{2}(\kappa)+N_{3}\kappa^{-2}
\Big( E_{x}^{\alpha_{\cdot}\bbeta ^{\rho}(\alpha_{\cdot},x)}
\sup_{t\leq\theta}|x_{t}-y_{t}|^{2}\Big)^{1/2},
\end{equation}
where $N_{1},N_{2},N_{3}$ are independent of $\alpha_{\cdot},\rho $, and $x$,
$N_{1}(\gamma)\to0$ as $\gamma\to\infty$,
  $N_{2}(\kappa)\to 0$ as $\kappa\downarrow0$, and we use the notation
$$
I^{h}(\alpha_{\cdot},\rho, x)
=E_{x}^{\alpha_{\cdot}\bbeta ^{\rho}(\alpha_{\cdot},x)}
\int_{0}^{\theta}|
h  (\alpha_{t},\beta (\alpha_{t},y_{t}),x_{t})-
h  (\alpha_{t},\beta (\alpha_{t},y_{t}),y_{t})|\,dt.
$$
By Corollary \ref{corollary 2.26.1}  the factor of $\gamma$
in \eqref{2.10.01}
is dominated by $\mu(\rho)$ for an appropriate function
$\mu(\rho)$ which tends to zero as $\rho\downarrow0$. The last term in \eqref{2.10.01}
is dominated by $\mu(\rho)\kappa^{-2}$.

After that taking into account \eqref{2.10.4} and Theorem 2.2.2 of \cite{Kr77}
we see that
$$
E_{x}^{\alpha_{\cdot}\bbeta ^{\rho}(\alpha_{\cdot},x)}
u_{K}(x_{\theta})e^{-\phi_{\theta}}\leq u_{K}(x)+\mu(\rho)(
 \gamma+\kappa^{-2})
+N_{1}(\gamma)+N_{2}(\kappa) 
$$
\begin{equation}
                                                              \label{2.10.5}
-
E_{x}^{\alpha_{\cdot}\bbeta ^{\rho}(\alpha_{\cdot},x)}\int_{0}^{\theta}
f(\alpha_{t},\beta (\alpha_{t},y_{t}),y_{t}) e^{-\phi_{t}}\,dt
+N m^{-1},
\end{equation}
where $N$ depend only on $d,\nu,K_{0}$ and the diameter of $G$.
We can replace the last $y_{t}$ in the integrand in \eqref{2.10.5}
by $x_{t}$ incurring as in
Corollary \ref{corollary 2.26.1} another error term like $\mu(\rho)  $
which goes to zero as $\rho\downarrow0$.
By adding to this that
$$
E_{x}^{\alpha_{\cdot}\bbeta ^{\rho}(\alpha_{\cdot},x)}
g(x_{\tau})e^{-\phi_{\tau}}=
E_{x}^{\alpha_{\cdot}\bbeta ^{\rho}(\alpha_{\cdot},x)}
u_{K}(x_{\tau})e^{-\phi_{\tau}}
$$
$$
\leq E_{x}^{\alpha_{\cdot}\bbeta ^{\rho}(\alpha_{\cdot},x)}
u_{K}(x_{\theta})e^{-\phi_{\theta}}+\sup_{G}|Du_{K}|
E_{x}^{\alpha_{\cdot}\bbeta ^{\rho}(\alpha_{\cdot},x)}
 |x_{\tau}-x_{\theta}|
$$
$$
+
\sup_{G}|u_{K}|
E_{x}^{\alpha_{\cdot}\bbeta ^{\rho}(\alpha_{\cdot},x)}
(\tau-\theta)
\leq E_{x}^{\alpha_{\cdot}\bbeta ^{\rho}(\alpha_{\cdot},x)}
u_{K}(x_{\theta})e^{-\phi_{\theta}}
$$
$$
+N_{4}
E_{x}^{\alpha_{\cdot}\bbeta ^{\rho}(\alpha_{\cdot},x)}
( (\tau -\theta)^{1/2}+\tau -\theta),
$$
where $N_{4}$ depends only on $u_{K}$, $d$, and $K_{0}$, 
we see that to prove \eqref{2.10.6}
it suffices now to show that
$$
\chi(\rho):=\sup_{\alpha_{\cdot}\in\frA}
E_{x}^{\alpha_{\cdot}\bbeta ^{\rho}(\alpha_{\cdot},x)} \int_{\theta}^{\tau}
 \,dt\to0
$$ 
as $\rho\downarrow 0$
uniformly with respect to $x$. By Lemma \ref{lemma 2.25.1}
 and It\^o's formula we have
$$
\chi(\rho)\leq
E_{x}^{\alpha_{\cdot}\bbeta ^{\rho}(\alpha_{\cdot},x)}
 \big[\Phi(x_{\theta})-\Phi(x_{\tau})\big] 
=E_{x}^{\alpha_{\cdot}\bbeta ^{\rho}(\alpha_{\cdot},x)}
\big[\Phi(x_{\theta})-\Phi(y_{\theta})\big]I_{\theta<\tau}
$$
and it only remains to use Lemma \ref{lemma 2.9.2}
once more. The theorem is proved.

{\bf Proof of Theorem \ref{theorem 2.10.1}}. First choose and fix
$K$ and $m$ so that $|v-u_{K}|\leq\varepsilon/4$ and $Nm^{-1}
\leq\varepsilon/4$, where $N$ is taken from   Theorem \ref{theorem 1.11.1}.
Then   find 
and fix $\kappa$ and $\gamma$  from $N_{1}(\gamma)
+N_{2}(\kappa)\leq \varepsilon/4$. Finally find
$\rho $ such that
$$
\mu(\rho)(1+\gamma+\kappa^{-2})\leq\varepsilon/4.
$$                            
  Then
\eqref{2.10.6} will become \eqref{2.10.8}. 

The last statement of the theorem follows by construction of $\bbeta^{\rho} 
(\alpha_{\cdot},x)$. The theorem is proved.

\begin{remark}
                                               \label{remark 2.16.1}
An important particular case of Theorem \ref{theorem 2.10.1}
is when $\sigma,b,c,f$ are independent of $\alpha$, so that we are actually dealing with
a controlled diffusion process. Also, clearly, similar statements to
Theorem \ref{theorem 2.10.1} hold true if we exchange the roles of $\alpha $ and 
 $\beta $ and consider the stochastic differential game corresponding to 
$$
H[u](x)=\infsup_{\beta\in B\,\,\alpha\in A}
[L^{\alpha\beta} u(x)+f^{\alpha\beta} (x)],
$$
in place of \eqref{1.16.1}. Of course, one should then replace
Assumption \ref{assumption 2.6.1} with a similar one about $B$.
To reduce this game to the one
we are treating, it suffices just to rename $A$ and $B$ and take $-u$, $-g$ and $-f$
in place of $u$, $g$, and $f$, respectively. 

\end{remark}

{\bf Proof of Theorem \ref{theorem 2.12.2}}. Fix $\nu>0$
and replace \eqref{5.11.1}
with 
$$
x_{t}=x+\int_{0}^{t}\sigma^{\alpha_{s}
\beta_{s} }(  x_{s})\,dw_{s}+\nu \hat w_{t}
+\int_{0}^{t}b^{\alpha_{s}
\beta_{s} }(  x_{s})\,ds.
$$
The solution of this equation is denoted by 
$x_{t}^{\alpha_{\cdot}\beta_{\cdot}x}(\nu)$
 and by $\tau^{\alpha_{\cdot}\beta_{\cdot}x}(\nu)$ we denote
its first exit time from $G$. We take the same $c,f,g$ and define
$v(x,\nu)$ by \eqref{2.12.2} where we replace $x_{t},\tau$, and $\phi_{t}$
with $x_{t}(\nu),\tau(\nu)$, and
$$
\phi_{t}(\nu)=\int_{0}^{t}c^{\alpha_{s}\beta_{s}}(x_{s}(\nu))\,ds,
$$
respectively. Obviously to thus obtained new stochastic differential game
we can apply Theorem \ref{theorem 2.10.1} and conclude that
for any $\varepsilon>0$ there exists $\beta(\alpha,x)$, with the properties
described in Theorem \ref{theorem 2.10.1} and $\rho_{0}>0$ such that
 if for $\rho\in(0,\rho_{0}]$,
 $x\in G$, and $\alpha_{\cdot}\in\frA$
we define the process $y_{t}=y_{t}^{\alpha_{\cdot}x}(\rho)$
as a solution of
\begin{equation}
                                                            \label{2.27.1}
dy_{t}=\sigma ^{(\rho)}(\alpha_{t},y_{t},y_{t})\,dw_{t}+\nu \hat w_{t}
+b ^{(\rho)}(\alpha_{t},y_{t},y_{t})\,dt,\quad t\geq0,\quad y_{0}=x,
\end{equation}
 then
$$
\sup_{\alpha_{\cdot}\in\frA} E_{x}^{\alpha_{\cdot}\bbeta^{\rho}  (\alpha_{\cdot},x)}\big[\int_{0}^{\tau(\nu)}
 f( x_{t}(\nu))e^{-\phi_{t}(\nu)}\,dt
$$
\begin{equation}
                                                              \label{2.27.2}  
+g(x_{\tau(\nu)}(\nu))e^{-\phi_{\tau(\nu)}(\nu)
 }\big]\leq v(x,\nu)+\varepsilon.
\end{equation}

It follows that to prove the theorem it suffices to show that
$$
E^{\alpha_{\cdot}\beta_{\cdot}}_{x}\Big[\int_{0}^{\tau}
f(x_{t})e^{-\phi_{t}}\,dt+g(x_{\tau})e^{-\phi_{\tau}}\Big]
$$
\begin{equation}
                                                      \label{2.27.3}
-E^{\alpha_{\cdot}\beta_{\cdot}}_{x}\Big[\int_{0}^{\tau(\nu)}
f(x_{t}(\nu))e^{-\phi_{t}(\nu)}\,dt+g(x_{\tau(\nu)(\nu)})
e^{-\phi_{\tau(\nu)}(\nu)}\Big]\to0
\end{equation}
as $\nu\downarrow0$ uniformly with respect to $\alpha_{\cdot}
\in\frA$, $\beta_{\cdot}\in\frB$, and $x\in G$.

First observe (although this is an overkill) that Lemma \ref{lemma 1.17.1} is applicable
here when $\sigma^{i}$'s are independent of the first 
space variable. Then  Corollary \ref{corollary 2.17.1}
is also applicable which as in Lemma \ref{lemma 2.9.2}
leads to the conclusion that
\begin{equation}
                                                      \label{2.27.03}
 E^{\alpha_{\cdot}\beta_{\cdot}}_{x}\sup_{t\leq\theta(\nu)}
|x_{t}-x_{t}(\nu)|^{2}\to0 
\end{equation}
as $\nu\downarrow0$ uniformly with respect to $\alpha_{\cdot}
\in\frA$, $\beta_{\cdot}\in\frB$, and $x\in G$,
where $\theta(\nu)$ is the minimum of exit times of $x_{t}$
and $x_{t}(\nu)$ from $G$.

Next, while proving \eqref{2.27.3} first assume that $g\equiv0$.
Observe that, owing to \eqref{2.27.03}, the argument at the end of the proof of
Theorem \ref{theorem 1.11.1} shows that it suffices to prove the version of
\eqref{2.27.3} when both $\tau$ and $\tau(\nu)$ are replaced with
$\theta(\nu)$ (assuming $g\equiv0$).

Then notice that in light of the continuity of $f$ in $x$
uniform with respect to $(\alpha,\beta)$ (cf.~also \eqref{2.27.4})
$$
E^{\alpha_{\cdot}\beta_{\cdot}}_{x} \int_{0}^{\theta(\nu)}
|f(x_{t})-f(x_{t}(\nu)|e^{-\phi_{t}} \,dt
\leq E^{\alpha_{\cdot}\beta_{\cdot}}_{x} \int_{0}^{\theta(\nu)}
|f(x_{t})-f(x_{t}(\nu)| \,dt\to 0
$$
as $\nu\downarrow0$ uniformly with respect to $\alpha_{\cdot}
\in\frA$, $\beta_{\cdot}\in\frB$, and $x\in G$. 

Also
 $$
I_{x}^{\alpha_{\cdot}\beta_{\cdot}}:=
E_{x}^{\alpha_{\cdot}\beta_{\cdot}} \int_{0}^{\theta(\nu)}
| f ( x_{t}(\nu))|\,|e^{-\phi_{t}}-e^{-\phi_{t}(\nu)}|\,dt 
$$ 
$$
\leq K_{0}E_{x}^{\alpha_{\cdot}\beta_{\cdot}}  \int_{0}^{\theta(\nu)}
\int_{0}^{t}| c ( x_{s})-
c ( x_{s}(\nu))|\,ds\,dt
$$
$$
=K_{0}E_{x}^{\alpha_{\cdot}\beta_{\cdot}} 
 \int_{0}^{\theta(\nu)}(\theta(\nu)-s)
| c ( x_{s})-
c ( x_{s}(\nu))|\,ds 
$$
$$
\leq K_{0}\Big[E_{x}^{\alpha_{\cdot}\beta_{\cdot}}\theta^{3}(\nu)\Big]^{1/2}
\Big[J_{x}^{\alpha_{\cdot}\beta_{\cdot}}\Big]^{1/2},
$$
where
$$
J_{x}^{\alpha_{\cdot}\beta_{\cdot}}=
 E_{x}^{\alpha_{\cdot}\beta_{\cdot}}
 \int_{0}^{\theta(\nu)}
 | c ( x_{s})-
c ( x_{s}(\nu))|^{2}\,ds \leq 2K_{0}
E_{x}^{\alpha_{\cdot}\beta_{\cdot}}
 \int_{0}^{\theta(\nu)}
 | c ( x_{s})-
c ( x_{s}(\nu))| \,ds.
$$
One sees easily as above that 
$I_{x}^{\alpha_{\cdot}\beta_{\cdot}}\to0$
as  $\nu\downarrow 0$ uniformly
with respect  to $\alpha_{\cdot}
\in\frA$, $\beta_{\cdot}\in\frB$, and $x\in G$.

 It remains to deal with the terms containing $g$ in \eqref{2.27.3}.
Since $g\in C^{2}(\bar G)$, by It\^o's formula we have
$$
E_{x}^{\alpha_{\cdot}\beta_{\cdot}} g(x_{\tau})e^{-\phi_{\tau}}
=g(x)
$$
\begin{equation}
                                                   \label{2.28.1}
+E_{x}^{\alpha_{\cdot}\beta_{\cdot}}\int_{0}^{\tau}
\big[a _{ij}(x_{t})D_{ij}g(x_{t})
+b _{i}(x_{t})D_{i}g(x_{t})-
c (x_{t})g(x_{t})\big]e^{-\phi_{t}}\,dt,
\end{equation}
$$
E_{x}^{\alpha_{\cdot}\beta_{\cdot}}g(x_{\tau(\nu)(\nu)})
e^{-\phi_{\tau(\nu)}(\nu)}
=g(x)+(1/2)\nu^{2}
E_{x}^{\alpha_{\cdot}\beta_{\cdot}}\int_{0}^{\tau(\nu)}
\Delta g(x_{t}(\nu)))e^{-\phi_{t}(\nu)}\,dt
$$
$$
+E_{x}^{\alpha_{\cdot}\beta_{\cdot}}\int_{0}^{\tau(\nu)}
\big[a _{ij}(x_{t}(\nu))D_{ij}g(x_{t}(\nu))
+b _{i}(x_{t}(\nu))D_{i}g(x_{t}(\nu))
$$
\begin{equation}
                                                   \label{2.28.2}
-
c (x_{t}(\nu))g(x_{t}(\nu))\big]e^{-\phi_{t}(\nu)}\,dt.
\end{equation}
The second term  on the right in \eqref{2.28.2}
clearly goes to zero as  $\nu\downarrow 0$ uniformly
with respect  to $\alpha_{\cdot}
\in\frA$, $\beta_{\cdot}\in\frB$, and $x\in G$.
The difference of the remaining ones in \eqref{2.28.1}
and \eqref{2.28.2} is shown to do the same by the first
part of the proof.  The theorem is proved.

{\bf Proof of Lemma \ref{lemma 2.21.1}}.
This proof if very similar to the second part of the proof
of Theorem \ref{theorem 2.12.2}.
 First we assume that $g\equiv 0$. Take
 $\theta=\theta^{\alpha_{\cdot}\bbeta ^{\rho}(\alpha_{\cdot},x)x}(\rho)$
from Lemma \ref{lemma 2.9.2} and note that the argument at the end of the proof of
Theorem \ref{theorem 1.11.1} shows that it suffices to prove the version of
\eqref{2.21.1} when both $\tau$ and $\tau(\rho)$ are replaced with
$\theta(\rho)$ (assuming $g\equiv0$).

Next,  observe that
 $$
E_{x}^{\alpha_{\cdot}\bbeta^{(\rho)}(\alpha_{\cdot},x)} \int_{0}^{\theta(\rho)}
| f( x_{t})-  f^{(\rho)} (\alpha_{t}, y_{t}(\rho))|e^{-\phi_{t}}\,dt 
$$ 
\begin{equation}
                                                           \label{2.22.2}
\leq  E_{x}^{\alpha_{\cdot}\bbeta^{(\rho)}(\alpha_{\cdot},x)} \int_{0}^{\theta(\rho)}
| f(\alpha_{t},\beta(\alpha_{t},y_{t}(\rho)), x_{t})-
f^{(\rho)} (\alpha_{t}, y_{t}(\rho), y_{t}(\rho))| \,dt .
\end{equation}
By Corollary \ref{corollary 2.26.1}
the last expression tends to zero
as $\rho\downarrow 0$ uniformly
with respect to $\alpha_{\cdot}\in \frA$ and $x\in G$.

Also as in the above proof
 $$
I_{x}^{\alpha_{\cdot}\bbeta^{(\rho)}(\alpha_{\cdot},x)}:=
E_{x}^{\alpha_{\cdot}\bbeta^{(\rho)}(\alpha_{\cdot},x)} \int_{0}^{\theta(\rho)}
|f^{(\rho)} (\alpha_{t}, y_{t}(\rho))|\,|e^{-\phi_{t}}-e^{-\phi_{t}(\rho)}|\,dt 
$$ 
$$
\leq K_{0}\Big[E_{x}^{\alpha_{\cdot}\bbeta^{(\rho)}(\alpha_{\cdot},x)}\theta^{3}(\rho)\Big]^{1/2}
\Big[J_{x}^{\alpha_{\cdot}\bbeta^{(\rho)}(\alpha_{\cdot},x)}\Big]^{1/2},
$$
where $J_{x}^{\alpha_{\cdot}\bbeta^{(\rho)}(\alpha_{\cdot},x)}$
stands for
$$
E_{x}^{\alpha_{\cdot}\bbeta^{(\rho)}(\alpha_{\cdot},x)}
 \int_{0}^{\theta(\rho)}
 | c(\alpha_{s},\beta(\alpha_{s},y_{s}(\rho)), x_{s})-
c^{(\rho)}(\alpha_{s},  y_{s})|^{2}\,ds 
$$
$$
\leq 2K_{0}E_{x}^{\alpha_{\cdot}\bbeta^{(\rho)}(\alpha_{\cdot},x)}
 \int_{0}^{\theta(\rho)}
 | c(\alpha_{s},\beta(\alpha_{s},y_{s}(\rho)), x_{s})-
c^{(\rho)}(\alpha_{s},   y_{s})| \,ds .
$$
Lemma \ref{lemma 2.25.2}
and Corollary \ref{corollary 2.26.1}
convince us that $I_{x}^{\alpha_{\cdot}\bbeta^{(\rho)}(\alpha_{\cdot},x)}\to0$
as  $\rho\downarrow 0$ uniformly
with respect to $\alpha_{\cdot}\in \frA$ and $x\in G$.

It remains to deal with the terms containing $g$ in \eqref{2.21.1}.
Again by using   It\^o's formula
 we write
$$
E_{x}^{\alpha_{\cdot}\bbeta^{(\rho)}(\alpha_{\cdot},x)} g(x_{\tau})e^{-\phi_{\tau}}
=g(x)
$$
$$
+E_{x}^{\alpha_{\cdot}\bbeta^{(\rho)}(\alpha_{\cdot},x)}\int_{0}^{\tau}
\big[a_{ij}(x_{t})D_{ij}g(x_{t})+b_{i}(x_{t})D_{i}g(x_{t})-c(x_{t})g(x_{t})\big]e^{-\phi_{t}}\,dt.
$$
Similarly we transform the term with $g$ involving $\tau(\rho)$
and then we reduce the problem to estimating the terms like
the ones we started with. The lemma is proved.

\mysection{A particular case where $A$ is a singleton}
                                                      \label{section 2.23.1}
Here we assume that $A$ is a singleton and will not write
$\alpha$ and $\alpha_{\cdot}$ in our notation. In particular, now we are dealing
with a controlled diffusion process given as a solution of the
equation
\begin{equation}
                                                                  \label{2.23.1}
dy_{t}=\sigma(\beta_{t},y_{t})\,dw_{t}+b(\beta_{t},y_{t})\,dt,\quad t\geq0,\quad
y_{0}=x.
\end{equation}
Its solution is denoted by $y^{\beta_{\cdot}x}_{t}$. Our goal is
to minimize
\begin{equation}
                                                          \label{2.23.6}
E^{\beta_{\cdot}}_{x}\Big[
\int_{0}^{\tau}f(y_{t})e^{-\phi_{t}}\,dt+g(y_{\tau})
e^{-\phi_{\tau}}\Big]
\end{equation}
over $\beta_{\cdot}\in\frB$, where (according to our standard notation)
$\tau^{\beta_{\cdot}x}$ is the first exit time of $y^{\beta_{\cdot}x}_{t}$
from $G$, $f(y_{t})=f(\beta_{t},y^{\beta_{\cdot}x}_{t})$,
$$
\phi^{\beta_{\cdot}x}_{t}=\int_{0}^{t}
c(\beta_{s},y^{\beta_{\cdot}x}_{s})\,ds.
$$

In this case Theorem \ref{theorem 2.10.1} becomes the following.

\begin{theorem}
                                       \label{theorem 2.23.1}
Under the assumptions of Theorem \ref{theorem 2.10.1}
for any $\varepsilon>0$ there exist
a Borel measurable
$B$-valued function $\beta( x)$ on $  \bR^{d}$
and $\rho_{0}>0$ such that, if for $\rho\in(0,\rho_{0}]$,
we define
$$
\sigma^{(\rho)}(z,y)=\int_{\bR^{d}}\sigma( \beta ( z+
\rho \xi),y+\rho \xi)\zeta(z)\,d\xi,
$$
introduce $b^{(\rho)}(z,y)$ similarly, and for
 $x\in G$ 
  define the process $z_{t}=z_{t}^{ x}(\rho)$ by
\begin{equation}
                                                       \label{2.28.5}
dz_{t}=\sigma^{(\rho)}(z_{t},z_{t})\,dw_{t}+
b^{(\rho)}(z_{t},z_{t})\,dt,\quad t\geq0,\quad z_{0}=x,
\end{equation}
and set $\beta^{\rho}_{t}( x)=\beta(z^{x}_{t}(\rho) )$, then
$$
\inf_{\beta_{\cdot}\in\frB}
E_{x}^{\beta_{\cdot} }\big[\int_{0}^{\tau}
 f( y_{t})e^{-\phi_{t}
 }\,dt+g(y_{\tau})e^{-\phi_{\tau}
 }\big]
$$
\begin{equation}
                                                       \label{2.23.2}
\geq 
E_{x}^{\beta^{\rho}_{\cdot}(x) }\big[\int_{0}^{\tau}
 f( y_{t})e^{-\phi_{t}
 }\,dt+g(y_{\tau})e^{-\phi_{\tau}
 }\big] -\varepsilon.
\end{equation}
\end{theorem} 

Here is a version of Theorem \ref{theorem 2.12.2}
\begin{theorem}
                                                    \label{theorem 2.28.1}
In Theorem \ref{theorem 2.23.1} drop Assumption \ref{assumption 2.12.1} but
suppose that on $(\Omega,\cF,P)$ there is
a Wiener process $(\hat w_{t}, \cF_{t}), t\geq0$,
independent of $w_{t}$.
 Then for any $\varepsilon>0$ there exists a constant $\nu>0$ 
such that all assertions of Theorem \ref{theorem 2.10.1} hold true
if we add to the right-hand side 
of \eqref{2.28.5} the term $\nu\, d\hat w_{t}$.

\end{theorem}

\begin{remark}
                                        \label{remark 2.23.1}
In Section \ref{section 2.28.3} we are going to maximize
\eqref{2.23.6} instead of minimizing it. One problem
is reduced to another just by changing signs of $f$ and $g$.
Also it is worth noting that in Section \ref{section 2.28.3} the parameter
used in maximization is called $\alpha_{\cdot}$ instead
of $\beta_{\cdot}$.

\end{remark}

\mysection{Adjoint $\varepsilon$-optimal Markov policies
for $\alpha$}
                                                      \label{section 2.28.3}

Take $\varepsilon>0$, $\rho>0$, $\beta(\alpha,x)$ from
Theorem \ref{theorem 2.10.1} use the notation \eqref{2.9.5}
and, for $\alpha_{\cdot}\in\frA$ and $x\in\bR^{d}$,
 defined the controlled diffusion process
$y_{t}(\rho)=y_{t}^{\alpha_{\cdot}x}(\rho)$ by
\begin{equation}
                                                   \label{2.28.6}
dy_{t}=\sigma^{(\rho)}(\alpha_{t},y_{t})\,dw_{t}
+b^{(\rho)}(\alpha_{t},y_{t})\,dt,\quad t\geq0,\quad y_{0}=x,
\end{equation}
with the reward function
\begin{equation}
                                                   \label{2.28.06}
E^{\alpha_{\cdot}}_{x}\Big[\int_{0}^{\tau(\rho)}
f(y_{t}(\rho))e^{-\phi_{t}(\rho)}\,dt+
g(y_{\tau(\rho)}(\rho))e^{-\phi_{\tau(\rho)}(\rho)}\Big].
\end{equation}
We are going to maximize \eqref{2.28.06} treating $\alpha$ here as $\beta$
in Section \ref{section 2.23.1} and adjusting the maximization
problem to the one of minimization.

However, there is a formal objection to overcome before we can translate
the results of Section \ref{section 2.23.1} to our situation. Namely,
in Section \ref{section 2.23.1}, the functions
$\sigma,b,c,f$ as inherited from taking $A$ as a singleton
were assumed to be continuous with respect to $\beta$. Therefore,
here we need our $\sigma^{(\rho)},b^{(\rho)},c^{(\rho)},f^{(\rho)}$
to be continuous with respect to $\alpha$ and they may fail to
be such because, even if $h$ in \eqref{2.9.5} is continuous
in the first argument $\alpha$ uniformly with respect to $\beta$,
$\beta(\alpha,y+\rho z)$ can be discontinuous as a function of $\alpha$.
Indeed, for different $\alpha$, $\beta(\alpha,x)$ can be   very different
functions of $x$. However, in light of the second statement
in Theorem \ref{theorem 2.10.1} to make $\beta(\alpha,x)$
 continuous with respect to $\alpha$ it suffices just to change
the distance function in $A$ keeping it the same as $\alpha_{1},\alpha_{2}$
belong to the same $A_{i}$ and defining it as $1$ otherwise. 
By the way, this change in no way affects the set of policies
of $\alpha$ and only allows us to formally apply the results
of Section \ref{section 2.23.1}.

According to Theorem \ref{theorem 2.23.1} for any $\varepsilon>0$
there exist a Borel measurable $A$-valued function
$\alpha^{\varepsilon}(z)$ on $\bR^{d}$ and a Lipschitz continuous functions
$\hat \sigma(z)$ and $\hat b(z)$ on $\bR^{d}$ with values
in the set of $d\times d_{1}$-matrices and in $\bR^{d}$,
respectively, such that, if for $x\in G$ we define the process
$z^{x}_{t}$ by
\begin{equation}
                                                       \label{3.1.1}
dz_{t}=\hat\sigma (z_{t} )\,dw_{t}+
\hat b  ( z_{t})\,dt,\quad t\geq0,\quad z_{0}=x,
\end{equation}
and set $\alpha^{\varepsilon,x}_{t}( x)=\alpha^{\varepsilon}(z^{x}_{t} )$, then
$$
\sup_{\alpha_{\cdot}\in\frA}
E_{x}^{\alpha_{\cdot}}\Big[\int_{0}^{\tau(\rho)}
f(y_{t}(\rho))e^{-\phi_{t}(\rho)}\,dt+
g(y_{\tau(\rho)}(\rho))e^{-\phi_{\tau(\rho)}(\rho)}\Big]
$$
\begin{equation}
                                                       \label{3.1.2}
\leq 
E_{x}^{\alpha^{\varepsilon,x}_{\cdot} }\Big[\int_{0}^{\tau(\rho)}
f(y_{t}(\rho))e^{-\phi_{t}(\rho)}\,dt+
g(y_{\tau(\rho)}(\rho))e^{-\phi_{\tau(\rho)}(\rho)}\Big] +\varepsilon.
\end{equation}

Finally, due to  Lemma \ref{lemma 2.21.1}, \eqref{3.1.2} implies that
\eqref{3.1.3} holds with $3\varepsilon$ in place of
$\varepsilon$.
This proves part (a) of Theorem \ref{theorem 2.10.2}.
The proof of part (b) is quite similar and the theorem is proved.

\mysection{Proof of Theorem \protect\ref{theorem 3.4.1}}
                                                        \label{section 3.10.4}

If in Theorem 14.1.6 of \cite{Kr_18} we replace $H[u]$ 
and $P[u]$ by $-H[-u]$ and $-P[-u]$, then we will see
that for any $K>0$ the equation
$$
\min(H[u_{K}],-P[-u_{-K}]+K)=0
$$
in $G$ (a.e.) with boundary condition $u_{K}=g \in C^{2}$
has a solution $u_{-K}\in W^{2}_{p}(G)$ for any $p>1$. By following the arguments
in  Section 7 of
\cite{Kr_14}, we conclude that $u_{-K}\uparrow v$ uniformly on 
$\bar G$  as $K\to\infty$. Observe that (a.e.) in $G$
\begin{equation}
                                                             \label{3.2.1}
H[u_{-K}]\geq0.
\end{equation}

Fix $K>0$ and $m\in\{1,2,...\}$.
In the same way in which we found above the function
$\beta(x)$ we find a Borel $A$-valued function
$\alpha(x)$ such that in $G$
$$
\inf_{\beta\in B}
[L^{\alpha(\cdot)\beta}u_{-K}+f^{\alpha(\cdot)\beta}]\geq -1/m.
$$
Our goal is to prove that if $K$ and $m$ are large enough and $\rho$
is small enough, then the above $\alpha(x)$ is the one we
are talking about in Theorem \ref{theorem 3.4.1}.

 Take $y^{\beta_{\cdot}x}_{t}(\rho)$ and $\balpha^{\rho}_{t}(\beta_{\cdot},x)=\alpha( 
 y_{t}^{\beta_{\cdot}x}(\rho))$ from the statement of the theorem.
Introduce $\theta=\theta^{\balpha^{\rho} (\beta_{\cdot},x)
 \beta_{\cdot} x}(\rho)$ 
 as the minimum of the first exit times
of $x_{ t}^{\balpha^{\rho} (\beta_{\cdot},x)\beta_{\cdot} x}$
and of $y_{t} ^{\beta_{\cdot}x}(\rho)$  from $G$. Then
in the same way in which we arrived at Lemma \ref{lemma 2.9.2} we obtain that
$$
\sup_{\beta_{\cdot}\in\frB}E_{x}^{\balpha^{\rho} (\beta_{\cdot},x)\beta_{\cdot}}
 \sup_{t\leq \theta(\rho)}|x_{t}-y_{t}(\rho)|^{2} \to0
$$
as $\rho\downarrow0$ uniformly with respect to $x\in G$.

Then following closely the argument in Section \ref{section 3.10.1}
we get an analog of Theorem \ref{theorem 1.11.1} that
for any
$x\in G$,  $\rho,\gamma,\kappa >0$   we have
$$
u_{-K}(x)\leq \inf_{\beta_{\cdot}\in\frB}
E_{x}^{\balpha^{\rho}(\beta_{\cdot},x)\beta_{\cdot}}\big[\int_{0}^{\tau}
 f( x_{t})e^{-\phi_{t}
 }\,dt+g(x_{\tau})e^{-\phi_{\tau}
 }\big] 
$$
$$
+\mu(\rho)(1+\gamma+\kappa^{-2}) +N_{1}(\gamma)+N_{2}(\kappa)+N m^{-1},
$$
where $N_{1}(\gamma)$ is independent of $\rho,\kappa $, 
$N_{1}(\gamma)\to0$ as $\gamma\to\infty$, $N_{2}(\kappa)$
is independent of $\rho$,
$N_{2}(\kappa)\to 0$ as $\kappa\downarrow0$, 
$N $ depends only on $d,\delta,K_{0}$, and the diameter of $G$,
$\mu(\rho)$ is independent of $\gamma,\kappa$
and $\mu(\rho)\to0$ as $\rho\downarrow0$.

After that the assertion of Theorem \ref{theorem 3.4.1}
is obtained by the same short argument as in Section
\ref{section 3.10.1} in the proof of Theorem \ref{theorem 2.10.1}.

\end{document}